\documentclass{article}
\usepackage{amsthm}
\usepackage{fleqn}
\usepackage{array}
\usepackage[latin1]{inputenc}
\usepackage{amssymb}
\usepackage{amsfonts}
\usepackage{amsmath}

\begin{document}

\newtheorem{thm}{Theorem}[section]
\newtheorem{prop}{Proposition}[section]
\newtheorem{lem}{Lemma}[section]
\newtheorem{cor}{Corollary}[section]

\numberwithin{equation}{section}

\begin{center}
{\Large Some classes of minimal surfaces 

\vspace{1mm}

in the $3$-space with $2m$-norm}

\vspace{5mm}

Makoto SAKAKI and Ryota TANAKA
\end{center}

\vspace{2mm}

{\bf Abstract.} We discuss translation minimal surfaces, homothetical minimal surfaces, and separable minimal surfaces in the $3$-space with $2m$-norm. 

\vspace{2mm}


\vspace{2mm}


\section{Introduction}

It is interesting to generalize differential geometry of surfaces in the Euclidean $3$-space to that in normed $3$-spaces (cf. [1], [2], [3], [4], [10], [11]), where the notion of Birkhoff orthogonality plays an important role. As special classes of minimal surfaces in the Euclidean $3$-space, translation and homothetical minimal surfaces are studed (cf. [8], [9]), and more generally, separable minimal surfaces are studied (cf. [5], [6], [7], [12]). In this paper, generalizing those results in the Euclidean $3$-space, we discuss translation minimal surfaces, homothetical minimal surfaces, and separable minimal surfaces in the $3$-space with $2m$-norm. 

In Section 2, following [2], we give basic facts on surfaces in normed $3$-spaces. In Section 3, we discuss graph surfaces in the $3$-space $({\mathbb R}^3, \|\cdot\|_{2m})$ with $2m$-norm. In Section 4, translation minimal surfaces in $({\mathbb R}^3, \|\cdot\|_{2m})$ are treated. In Section 5, homothetical minimal surfaces in $({\mathbb R}^3, \|\cdot\|_{2m})$ are considered. In Section 6, we discuss the construction and examples of separable minimal surfaces in $({\mathbb R}^3, \|\cdot\|_{2m})$.

\section{Surfaces in normed $3$-spaces}

Let $({\mathbb R}^3, \|\cdot\|)$ be a normed $3$-space whose unit ball $B$ and unit sphere $S$ are given by
\[B = \{x \in {\mathbb R}^3; \| x\| \leq 1\}, \ \ \ \ S = \{x \in {\mathbb R}^3; \| x\| = 1\}. \]
We assume that $S$ is smooth and strictly convex, that is, $S$ is a smooth surface and $S$ contains no line segment. 

Let $v$ be a non-zero vector in ${\mathbb R}^3$ and $\Pi$ be a plane in ${\mathbb R}^3$. We say that $v$ is Birkhoff orthogonal to $\Pi$, denoted by $v \dashv_{B} \Pi$, if the tangent plane of $S$ at $v/\|v\|$ is parallel to $\Pi$. 

Let $M$ be a surface immersed in $({\mathbb R}^3, \|\cdot\|)$. Let $T_{p}M$ denote the tangent plane of $M$ at $p \in M$. There exists a vector $\eta(p) \in S$ such that $\eta(p) \dashv_{B} T_{p}M$. This gives a local smooth map $\eta: U \subset M \rightarrow S$, that is called the Birkhoff-Gauss map. It can be global if and only if $M$ is orientable. 

The mean curvature $H$ of $M$ at $p$ is defined by
\[H = \frac{1}{2}\mbox{trace}(d\eta_p). \]
The surface $M$ is called minimal if $H = 0$ identically.

\section{Graph surfaces in the $3$-space with $2m$-norm}

In the following, let $({\mathbb R}^3, \|\cdot\|_{2m})$ be the $3$-space with $2m$-norm
\[\|x\|_{2m} = (x_{1}^{2m}+x_{2}^{2m}+x_3^{2m})^{\frac{1}{2m}}, \ \ \ x = (x_1, x_2, x_3) \in {\mathbb R}^3, \]
where $m$ is a positive integer. Set
\[\Phi(x_1, x_2, x_3) := x_{1}^{2m}+x_{2}^{2m}+x_3^{2m}. \]
Then the unit sphere $S$ is given by
\[S = \{(x_1, x_2, x_3) \in {\mathbb R}^3; \Phi(x_1, x_2, x_3) = 1 \}. \]
Since the case $m = 1$ is the Euclidean case, we assume that $m \geq 2$ in the following. 

Let $M$ be a graph surface in $({\mathbb R}^3, \|\cdot\|_{2m})$ which is given by
\[X(u, v) = (u, v, f(u, v)) \]
for a smooth function $f(u, v)$ with $f_u \neq0$ and $f_v \neq 0$. Then 
\[X_u = (1, 0, f_u), \ \ \ X_v = (0, 1, f_v). \]
The Birkhoff-Gauss map $\eta = \eta(u, v)$ is characterized by the condition 
\[(\mbox{grad}(\Phi))_{\eta} = \left( \frac{\partial \Phi}{\partial x_1}(\eta), \frac{\partial \Phi}{\partial x_2}(\eta), \frac{\partial \Phi}{\partial x_3}(\eta) \right) = \varphi X_u \times X_v, \]
where $\varphi$ is a positive function and $\times$ is the standard cross product in ${\mathbb R}^3$. Then we have
\[\eta = A^{-\frac{1}{2m}} ( -f_{u}^{\frac{1}{2m-1}}, -f_{v}^{\frac{1}{2m-1}}, 1 ) \]
where
\[A := 1+f_{u}^{\frac{2m}{2m-1}}+f_{v}^{\frac{2m}{2m-1}}. \]

We can compute that
\[\eta_u = -\frac{1}{2m-1}A^{-\frac{2m+1}{2m}} f_{u}^{-\frac{2m-2}{2m-1}} f_{v}^{-\frac{2m-2}{2m-1}} \times \hspace{4cm} \]
\[ \left[ \left\{ \left( f_{v}^{\frac{2m-2}{2m-1}}+f_{v}^{2} \right)f_{uu}-f_{u}f_{v}f_{uv} \right\}X_{u}+\left\{ \left( f_{u}^{\frac{2m-2}{2m-1}}+f_{u}^2 \right)f_{uv}-f_{u}f_{v}f_{uu} \right\}X_{v} \right] \]
and
\[\eta_v = -\frac{1}{2m-1}A^{-\frac{2m+1}{2m}} f_{u}^{-\frac{2m-2}{2m-1}} f_{v}^{-\frac{2m-2}{2m-1}} \times \hspace{4cm} \]
\[ \left[ \left\{ \left( f_{v}^{\frac{2m-2}{2m-1}}+f_{v}^{2} \right)f_{uv}-f_{u}f_{v}f_{vv} \right\}X_{u}+\left\{ \left( f_{u}^{\frac{2m-2}{2m-1}}+f_{u}^2 \right)f_{vv}-f_{u}f_{v}f_{uv} \right\}X_{v} \right]. \]
So the the mean curvature $H$ is given by
\begin{eqnarray}
H = -\frac{1}{2(2m-1)}A^{-\frac{2m+1}{2m}} f_{u}^{-\frac{2m-2}{2m-1}} f_{v}^{-\frac{2m-2}{2m-1}} \times \nonumber 
\end{eqnarray}
\begin{eqnarray}
\hspace{5mm} \left\{ \left( f_{v}^{\frac{2m-2}{2m-1}}+f_{v}^{2} \right)f_{uu}-2f_{u}f_{v}f_{uv}+\left( f_{u}^{\frac{2m-2}{2m-1}}+f_{u}^2 \right)f_{vv} \right\}. 
\end{eqnarray}

\vspace{2mm}

(i) When $f(u, v) = g(u)+h(v)$, the surface $M$ is called a translation surface, where we assmue that $g'(u) \neq 0$ and $h'(v) \neq 0$. In this case, by (3.1), we have
\begin{eqnarray}
H =  -\frac{1}{2(2m-1)} \left( 1+(g')^{\frac{2m}{2m-1}}+(h')^{\frac{2m}{2m-1}} \right)^{-\frac{2m+1}{2m}} (g')^{-\frac{2m-2}{2m-1}} (h')^{-\frac{2m-2}{2m-1}} \nonumber 
\end{eqnarray}
\begin{eqnarray}
\hspace{1cm} \times \left\{ \left( (h')^{\frac{2m-2}{2m-1}}+(h')^{2} \right) g''+\left( (g')^{\frac{2m-2}{2m-1}}+(g')^2 \right) h'' \right\}. 
\end{eqnarray}

\vspace{2mm}

(ii) When $f(u, v) = g(u)h(v)$, the surface $M$ is called a homothetical surface, where we assmue that
\[g(u) \neq 0, \ \ \ g'(u) \neq 0, \ \ \ h(v) \neq 0, \ \ \ h'(v) \neq 0. \]
In this case, by (3.1), we have
{\small 
\begin{eqnarray}
H = -\frac{1}{2(2m-1)} \left( 1+(g' h)^{\frac{2m}{2m-1}}+(gh')^{\frac{2m}{2m-1}} \right)^{-\frac{2m+1}{2m}} (g' h)^{-\frac{2m-2}{2m-1}} (gh')^{-\frac{2m-2}{2m-1}} \nonumber 
\end{eqnarray}
\begin{eqnarray}
\times \left\{ \left( (gh' )^{\frac{2m-2}{2m-1}}+(gh' )^{2} \right)g'' h-2gh(g')^{2}(h')^{2}+\left( (g' h)^{\frac{2m-2}{2m-1}}+(g' h)^{2} \right)gh'' \right\}. \hspace{3mm} 
\end{eqnarray} }

\section{Translation minimal surfaces}

Let $M$ be a translation surface in $({\mathbb R}^3, \|\cdot\|_{2m})$ as in Section 3. By (3.2), we see that $H = 0$ identically if and only if
\[\left( (h')^{\frac{2m-2}{2m-1}}+(h')^{2} \right) g''+\left( (g')^{\frac{2m-2}{2m-1}}+(g')^2 \right) h'' = 0. \]
Then we have
\[\frac{g''}{(g')^{\frac{2m-2}{2m-1}}+(g')^2} = -\frac{h''}{(h')^{\frac{2m-2}{2m-1}}+(h')^{2}} = a, \]
where $a$ is a constant. 

If $a = 0$, then $g$ and $h$ are both linear, so that $M$ is a plane. We assume that $a \neq 0$ in the following. By integration we have
\[\int \frac{ g''}{(g')^{\frac{2m-2}{2m-1}}+(g')^{2}} du = au, \ \ \ \int \frac{h''}{(h')^{\frac{2m-2}{2m-1}}+(h')^{2}} dv = -av, \]
where the integral constants are omitted since $a \neq 0$. Letting $(g')^{\frac{1}{2m-1}} =: s$ and $(h')^{\frac{1}{2m-1}} =: t$, we can get
\[au = \int \frac{ g''}{(g')^{\frac{2m-2}{2m-1}}+(g')^{2}} du = (2m-1) \int \frac{1}{1+s^{2m}}ds =: F_{m}(s) \]
and
\[-av = \int \frac{h''}{(h')^{\frac{2m-2}{2m-1}}+(h')^{2}} dv = (2m-1) \int \frac{1}{1+t^{2m}}dt = F_{m}(t). \]
Hence
\[F_{m}\left( (g')^{\frac{1}{2m-1}} \right) = au, \ \ \ F_{m}\left( (h')^{\frac{1}{2m-1}} \right) = -av, \]
and
\[g' = \{(F_m)^{-1}(au)\}^{2m-1}, \ \ \ \ h' = \{(F_m)^{-1}(-av)\}^{2m-1}. \]
Thus we have
\[g(u) = \int \{(F_m)^{-1}(au)\}^{2m-1} du \]
and
\[h(v) = \int \{(F_m)^{-1}(-av)\}^{2m-1} dv. \]
So the surface $M$ is parametrized as 
{\small 
\begin{eqnarray}
X(u, v) = \left( u, v, \int \{(F_m)^{-1}(au)\}^{2m-1} du+\int \{(F_m)^{-1}(-av)\}^{2m-1} dv \right). 
\end{eqnarray} }

Since
\[(F_m)^{-1}(au) = (g')^{\frac{1}{2m-1}} = s, \]
the integration by substitution implies that
\[g(u(s)) = \frac{2m-1}{a}\int\frac{s^{2m-1}}{1+s^{2m}}ds = \frac{2m-1}{2ma}\log{(1+s^{2m})}. \]
Analogously, by
\[(F_m)^{-1}(-av) = (h')^{\frac{1}{2m-1}} = t, \]
we have
\[h(v(t)) = -\frac{2m-1}{2ma}\log{(1+t^{2m})}. \]
Thus the surface $M$ is reparametrized as 
{\small 
\begin{eqnarray}
X(s, t) = \frac{1}{a}\left( F_{m}(s), \ -F_{m}(t), \ \frac{2m-1}{2m}\left\{ \log{(1+s^{2m})}-\log{(1+t^{2m})} \right\} \right). 
\end{eqnarray} }

Hence we get the following: 

\begin{thm}
Let $M$ be a translation surface in $({\mathbb R}^3, \|\cdot\|_{2m})$ as above. Then $M$ is minimal if and only if it is a plane or given by (4.1). The surface given by (4.1) is reparametrized as in (4.2). 
\end{thm}

{\bf Remark.} The equation (4.2) can be represented by elementary functions. For example, when $m = 2$, since
\[\frac{1}{1+s^{4}} = \frac{\sqrt{2}}{4}\left( \frac{s+\sqrt{2}}{s^{2}+\sqrt{2}s+1}-\frac{s-\sqrt{2}}{s^{2}-\sqrt{2}s+1} \right), \]
we have
\[F_{2}(s) = 3\int \frac{1}{1+s^{4}}ds  = \frac{3\sqrt{2}}{8}\left\{ \log{\left( s^{2}+\sqrt{2}s+1 \right)}-\log{\left( s^{2}-\sqrt{2}s+1 \right)} \right\} \]
\[+\frac{3\sqrt{2}}{4}\left\{ \arctan{\left( \sqrt{2}s+1 \right)}+\arctan{\left( \sqrt{2}s-1 \right)} \right\}. \]

\section{Homothetical minimal surfaces}

Let $M$ be a homothetical surface in $({\mathbb R}^3, \|\cdot\|_{2m})$ as in Section 3. By (3.3), we see that $H = 0$ identically if and only if 
\begin{eqnarray}
\hspace{5mm} \left( (gh' )^{\frac{2m-2}{2m-1}}+(gh' )^{2} \right)g'' h-2gh(g')^{2}(h')^{2} \nonumber 
\end{eqnarray}
\begin{eqnarray}
\hspace{2cm} +\left( (g' h)^{\frac{2m-2}{2m-1}}+(g' h)^{2} \right)gh'' = 0. 
\end{eqnarray}

\vspace{2mm}

(i) The case where $g'' = 0$. Then $g(u) = au+b$ for constants $a \neq 0$ and $b$. The equation (5.1) reduces to
\[-2a^{2}h(h')^{2}+h'' \left( (ah)^{\frac{2m-2}{2m-1}}+(ah)^{2} \right) = 0, \]
so that
\[\frac{h'' }{h' } = \frac{2a^{2}hh' }{a^{\frac{2m-2}{2m-1}}h^{\frac{2m-2}{2m-1}}+a^{2}h^{2}}. \]
By integration we have
\[\log{|h'|} = 2a^{2}\int\frac{h}{a^{\frac{2m-2}{2m-1}}h^{\frac{2m-2}{2m-1}}+a^{2}h^{2}}dh+c_{1} \]
where $c_1$ is a constant. Letting $h^{\frac{2}{2m-1}} =: t$, we have
\[\log{|h'|} = (2m-1)a^{2}\int\frac{t^{m-1}}{a^{\frac{2m-2}{2m-1}}+a^{2}t^{m}}dt+c_{1} \]
\[= \frac{2m-1}{m}\log{\left( a^{\frac{2m-2}{2m-1}}+a^{2}t^{m} \right)}+c_{1} \]
\[= \frac{2m-1}{m}\log{\left( a^{\frac{2m-2}{2m-1}}+a^{2}h^{\frac{2m}{2m-1}} \right)}+c_{1}, \]
and
\[\frac{dh}{dv} = c_{2}\left( a^{\frac{2m-2}{2m-1}}+a^{2}h^{\frac{2m}{2m-1}} \right)^{\frac{2m-1}{m}} \]
where $c_2$ is a non-zero constant. Then we have
\begin{eqnarray}
\hspace{5mm} v = \frac{1}{c_2}\int \frac{1}{\left( a^{\frac{2m-2}{2m-1}}+a^{2}h^{\frac{2m}{2m-1}} \right)^{\frac{2m-1}{m}}} dh =: \Psi(h). 
\end{eqnarray}
So the surface $M$ is parametrized as
\begin{eqnarray}
\hspace{5mm} X(u, v) = (u, v, (au+b)\Psi^{-1}(v)), 
\end{eqnarray}
and is reparametrized as
\begin{eqnarray}
\hspace{5mm} X(u, h) = (u, \Psi(h), (au+b)h). 
\end{eqnarray}
The case where $h'' = 0$ is analogous.

\vspace{2mm}

(ii) The case where $g'' \neq 0$ and $h'' \neq 0$. Dividing (5.1) by $gh(g')^{2}(h')^{2}$, we have
{\small 
\begin{eqnarray}
\frac{g'' }{g^{\frac{1}{2m-1}}(g')^{2}(h')^{\frac{2m}{2m-1}}}+\frac{gg''}{(g')^{2}}-2+\frac{h''}{(g')^{\frac{2m}{2m-1}}h^{\frac{1}{2m-1}}(h')^{2}}+\frac{hh''}{(h')^{2}} = 0. 
\end{eqnarray} }
Differentiating with respect to $u$ and then $v$, we have
\[\left( \frac{g'' }{g^{\frac{1}{2m-1}}(g')^{2}} \right)' \left( (h')^{-\frac{2m}{2m-1}} \right)'+\left( (g')^{-\frac{2m}{2m-1}} \right)' \left( \frac{h''}{h^{\frac{1}{2m-1}}(h')^{2}} \right)' = 0. \]
Noticing that $g'' h'' \neq 0$, we have
\[\left( \frac{g'' }{g^{\frac{1}{2m-1}}(g')^{2}} \right)' \frac{1}{\left( (g')^{-\frac{2m}{2m-1}} \right)'} = -\left( \frac{h''}{h^{\frac{1}{2m-1}}(h')^{2}} \right)' \frac{1}{\left( (h')^{-\frac{2m}{2m-1}} \right)'} = a, \]
where $a$ is a constant. Then by integration, 
\[\frac{g'' }{g^{\frac{1}{2m-1}}(g')^{2}} = \frac{a}{(g')^{\frac{2m}{2m-1}}}+b, \ \ \ \ -\frac{h''}{h^{\frac{1}{2m-1}}(h')^{2}} = \frac{a}{(h')^{\frac{2m}{2m-1}}}+c \]
for constants $b$ and $c$, equivalently, 
\begin{eqnarray}
\hspace{5mm} g'' = g^{\frac{1}{2m-1}} \left( a(g')^{\frac{2m-2}{2m-1}}+b(g')^{2} \right) 
\end{eqnarray}
and
\begin{eqnarray}
\hspace{5mm} h'' = -h^{\frac{1}{2m-1}}\left( a(h')^{\frac{2m-2}{2m-1}}+c(h')^{2} \right). 
\end{eqnarray}

Substituting (5.6) and (5.7) into (5.5), we obtain that 
\[\frac{c-ag^{\frac{2m}{2m-1}}} {(g')^{\frac{2m}{2m-1}}}+2-bg^{\frac{2m}{2m-1}} = \frac{b-ah^{\frac{2m}{2m-1}}} {(h')^{\frac{2m}{2m-1}}}-ch^{\frac{2m}{2m-1}} = \lambda, \]
where $\lambda$ is a constant, and then
\begin{eqnarray}
(g')^{\frac{2m}{2m-1}} = \frac{c-ag^{\frac{2m}{2m-1}}} {\lambda+bg^{\frac{2m}{2m-1}}-2}, \ \ \ \ (h')^{\frac{2m}{2m-1}} = \frac{b-ah^{\frac{2m}{2m-1}}} {\lambda+ch^{\frac{2m}{2m-1}}}. 
\end{eqnarray}
Differentiating with respect to $u$ and $v$, respectively, we can get
\begin{eqnarray}
\hspace{5mm} g'' = -g^{\frac{1}{2m-1}} (g')^{\frac{2m-2}{2m-1}} \frac{a(\lambda-2)+bc}{\left( \lambda+bg^{\frac{2m}{2m-1}}-2 \right)^{2}} 
\end{eqnarray}
and
\begin{eqnarray}
\hspace{5mm} h'' = -h^{\frac{1}{2m-1}} (h')^{\frac{2m-2}{2m-1}} \frac{a\lambda+bc}{\left( \lambda+ch^{\frac{2m}{2m-1}} \right)^{2}}. 
\end{eqnarray}

The right hand side of (5.6) is equal to that of (5.9), and using (5.8), we find that
\[(a(\lambda-2)+bc) \left( \lambda-1+bg^{\frac{2m}{2m-1}} \right) = 0. \]
Analogously, by (5.7), (5.10) and (5.8), we have
\[(a\lambda+bc)\left( \lambda-1+ch^{\frac{2m}{2m-1}} \right) = 0. \]
There are four cases: 

(a) If $a(\lambda-2)+bc = a\lambda+bc = 0$, then $a = 0$ and $bc = 0$, from which we have $g'' = h'' = 0$, a contradiction. 

(b) If $a(\lambda-2)+bc = 0$ and $a\lambda+bc \neq 0$, then $\lambda-1+ch^{\frac{2m}{2m-1}} = 0$, from which we have $c = \lambda-1 = a = 0$ and $g'' = h''= 0$, a contradiction. 

(c) If $a(\lambda-2)+bc \neq 0$ and $a\lambda+bc = 0$, then $\lambda-1+bg^{\frac{2m}{2m-1}} = 0$, from which we have $b = \lambda-1 = a = 0$ and $g'' = h'' = 0$, a contradiction. 

(d) If $a(\lambda-2)+bc \neq 0$ and $a\lambda+bc \neq 0$, then $\lambda-1+bg^{\frac{2m}{2m-1}} = 0$ and $\lambda-1+ch^{\frac{2m}{2m-1}} = 0$, from which we have $b = c = \lambda-1 = 0$. By (5.8), we have
\[(g')^{\frac{2m}{2m-1}} = ag^{\frac{2m}{2m-1}}, \ \ \ \ (h')^{\frac{2m}{2m-1}} = -ah^{\frac{2m}{2m-1}}, \]
so that $a = 0$ and $g' = h' = 0$, a contradiction. 

Therefore, the case (ii) does not occur, and we have shown the following: 

\begin{thm}
Let $M$ be a homothetical surface in $({\mathbb R}^3, \|\cdot\|_{2m})$ as above. Then $M$ is minimal if and only if it is given by (5.2) and (5.3), which is reparametrized as in (5.4). 
\end{thm}

{\bf Remark.} The discussion in this section is analogous to that in [9].

\section{Separable minimal surfaces}

Let $M$ be an implicit surface in $({\mathbb R}^3, \|\cdot\|_{2m})$ which is given by 
\[f(x_1)+g(x_2)+h(x_3) = 0 \]
for smooth functions $f(x_1)$, $g(x_2)$ and $h(x_3)$ with $(f'(x_1), g'(x_2), h'(x_3)) \neq (0, 0, 0)$. Such a surface is called a separable surface. Translation surfaces and homothetical surfaces are examples of separable surfaces. 

We assume that $M$ is not a translation surface. Then we may assume that $f'(x_1)g'(x_2)h'(x_3) \neq 0$ and $f''(x_1)g''(x_2)h''(x_3) \neq 0$. 

The Birkhoff-Gauss map $\eta$ is characterized by the condition 
\[(\mbox{grad}(\Phi))_{\eta} = \varphi (f', g', h'), \]
where $\varphi$ is a positive function. Then we have
\[\eta = A^{-\frac{1}{2m}}\left( (f')^{\frac{1}{2m-1}}, (g')^{\frac{1}{2m-1}}, (h')^{\frac{1}{2m-1}} \right) \]
where 
\[A := (f')^{\frac{2m}{2m-1}}+(g')^{\frac{2m}{2m-1}}+(h')^{\frac{2m}{2m-1}}. \]

Since $h'(x_3) \neq 0$, $M$ can be expressed as $x_3 = x_3(x_1, x_2)$ and
\[\frac{\partial x_3}{\partial x_1} = -\frac{f'}{h'}, \ \ \ \ \frac{\partial x_3}{\partial x_2} = -\frac{g'}{h'}. \]
So $M$ is parametrized as $p(x_1, x_2) = (x_1, x_2, x_3(x_1, x_2))$ and
\[p_{x_1} = \left( 1, 0, -\frac{f'}{h'} \right), \ \ \ \ p_{x_2} = \left(0, 1, -\frac{g'}{h'} \right). \]

With respect to these coordinates $x_1$, $x_2$, we can compute that $\eta_{x_1} = \eta_{1}^{1} p_{x_1}+\eta_{1}^{2} p_{x_2}$, where
\[\eta_{1}^{1} = \frac{A^{-\frac{2m+1}{2m}}} {2m-1} \left[ (f')^{\frac{2m}{2m-1}} (h')^{-\frac{2m-2}{2m-1}} h''+\left\{ (g')^{\frac{2m}{2m-1}}+(h')^{\frac{2m}{2m-1}} \right\} (f')^{-\frac{2m-2}{2m-1}} f'' \right] \]
and
\[\eta_{1}^{2} = \frac{A^{-\frac{2m+1}{2m}}} {2m-1} (g')^{\frac{1}{2m-1}} \left\{ f' (h')^{-\frac{2m-2}{2m-1}} h''-(f')^{\frac{1}{2m-1}} f'' \right\}. \]
Analogously, we have $\eta_{x_2} = \eta_{2}^{1} p_{x_1}+\eta_{2}^{2} p_{x_2}$, where
\[\eta_{2}^{1} = \frac{A^{-\frac{2m+1}{2m}}} {2m-1} (f')^{\frac{1}{2m-1}} \left\{ g' (h')^{-\frac{2m-2}{2m-1}} h''-(g')^{\frac{1}{2m-1}} g'' \right\} \]
and
\[\eta_{2}^{2} = \frac{A^{-\frac{2m+1}{2m}}} {2m-1} \left[ (g')^{\frac{2m}{2m-1}} (h')^{-\frac{2m-2}{2m-1}} h''+\left\{ (f')^{\frac{2m}{2m-1}}+(h')^{\frac{2m}{2m-1}} \right\} (g')^{-\frac{2m-2}{2m-1}} g''\right]. \]
Hence, $M$ is minimal if and only if 
\[\left\{ (g')^{\frac{2m}{2m-1}}+(h')^{\frac{2m}{2m-1}} \right\} (f')^{-\frac{2m-2}{2m-1}} f''+\left\{ (f')^{\frac{2m}{2m-1}}+(h')^{\frac{2m}{2m-1}} \right\} (g')^{-\frac{2m-2}{2m-1}} g'' \]
\begin{eqnarray}
\hspace{1cm} +\left\{ (f')^{\frac{2m}{2m-1}}+(g')^{\frac{2m}{2m-1}} \right\} (h')^{-\frac{2m-2}{2m-1}} h'' = 0. 
\end{eqnarray}

\vspace{1mm}

We will continue the discussion just as in [7]. Let us introduce new variables $u$, $v$ and $w$ as
\[u = f(x_1), \ \ \ v = g(x_2), \ \ \ w = h(x_3). \]
Set
\[X(u) = (f')^{\frac{2m}{2m-1}}, \ \ \ Y(v) = (g')^{\frac{2m}{2m-1}}, \ \ \ Z(w) = (h')^{\frac{2m}{2m-1}}. \]
Noticing that
\[X' = \frac{dX}{du} = \frac{dX}{dx_1} \frac{dx_1}{du} = \frac{2m}{2m-1}(f')^{-\frac{2m-2}{2m-1}} f'' \neq 0, \]
\[Y' = \frac{dY}{dv} = \frac{dY}{dx_2} \frac{dx_2}{dv} = \frac{2m}{2m-1}(g')^{-\frac{2m-2}{2m-1}} g'' \neq 0, \]
\[Z' = \frac{dZ}{dw} = \frac{dZ}{dx_3} \frac{dx_3}{dw} = \frac{2m}{2m-1}(h')^{-\frac{2m-2}{2m-1}} h'' \neq 0, \]
we can rewrite the equation (6.1) as 
\begin{eqnarray}
B(u, v, w) := (Y(v)+Z(w))X'(u)+(Z(w)+X(u))Y'(v) \nonumber 
\end{eqnarray}
\begin{eqnarray}
\hspace{3cm} +(X(u)+Y(v))Z'(w) = 0 
\end{eqnarray}
for all values $u$, $v$ and $w$ such that $u+v+w = 0$.

\vspace{2mm}

{\bf Remark.} The functions $X$, $Y$ and $Z$ depend on $m$, but the equation (6.2) itself is independent of $m$. 

\begin{lem}
(cf. [7]) \ Let $Q(u, v, w )$ be a smooth function on a domain $\Omega \subset {\mathbb R}^3$. Supposse that $Q(u, v, w) = 0$ on the section $\Omega\cap\Pi$ where $\Pi: u+v+w = 0$. Then $Q_u = Q_v = Q_w$ on $\Omega\cap\Pi$. 
\end{lem}

Applying the lemma to $B(u, v, w)$, we have $B_u-B_v = 0$, $B_v-B_w = 0$ and $B_u-B_w = 0$, which are written as 
\begin{eqnarray}
\hspace{5mm} L_1 := (Y+Z)X''-(Z+X)Y''+(X'-Y')Z' = 0, 
\end{eqnarray}
\begin{eqnarray}
\hspace{5mm} L_2 := (Y'-Z')X'+(Z+X)Y''-(X+Y)Z'' = 0, 
\end{eqnarray}
\begin{eqnarray}
\hspace{5mm} L_3 := (Y+Z)X''-(Z'-X')Y'-(X+Y)Z'' = 0. 
\end{eqnarray}
Here $L_3 = L_1+L_2$. From (6.2), (6.3) and (6.4), we have a system of linear equations on $Y+Z$, $Z+X$ and $X+Y$, whose determinant $D$ is given by
\[D = X'Y''Z''+X''Y'Z''+X''Y''Z'. \]
Solving the system, we get
\begin{eqnarray}
\hspace{5mm} D(Y+Z)X' = X'Y'Z' \left\{ Y''(Z'-X')-Z''(X'-Y') \right\}, 
\end{eqnarray}
\begin{eqnarray}
\hspace{5mm} D(Z+X)Y' = X'Y'Z' \left\{ Z''(X'-Y')-X''(Y'-Z') \right\}, 
\end{eqnarray}
\begin{eqnarray}
\hspace{5mm} D(X+Y)Z' = X'Y'Z' \left\{ X''(Y'-Z')-Y''(Z'-X') \right\}. 
\end{eqnarray}

Applying Lemma 6.1 again to $L_1$, $L_2$ and $L_3$, we have $(L_1)_v = (L_1)_w$, $(L_2)_u = (L_2)_w$ and $(L_3)_u = (L_3)_v$, which give
\[Z''(X'-Y')-X''(Y'-Z') = -Y^{(3)}(Z+X), \]
\[X''(Y'-Z')-Y''(Z'-X') = -Z^{(3)}(X+Y), \]
\[Y''(Z'-X')-Z''(X'-Y')= -X^{(3)}(Y+Z). \]
Substituting them into (6.6), (6.7) and (6.8), we have
\[D(Y+Z)X' = -X'Y'Z' (Y+Z)X^{(3)}, \]
\[D(Z+X)Y' = -X'Y'Z' (Z+X)Y^{(3)}, \]
\[D(X+Y)Z' = -X'Y'Z' (X+Y)Z^{(3)}. \]
Since $X'Y'Z' \neq 0$, we find that
\begin{eqnarray}
\hspace{1cm} \frac{X^{(3)}}{X'} = \frac{Y^{(3)}}{Y'} = \frac{Z^{(3)}}{Z'} = a 
\end{eqnarray}
where $a$ is a constant, because $X^{(3)}/X'$ is a function of $u$ and $Y^{(3)}/Y'$ is a function of $v$.

\vspace{2mm}

(i) The case where $a > 0$. Let $a = b^2$ for $b > 0$. The solutions of (6.9) are given by
\begin{eqnarray}
\hspace{5mm} \begin{array}{c}
X(u) = p_{1}+q_{1}e^{bu}+r_{1}e^{-bu}, \\ \vspace{1mm} 
Y(v) =p_{2}+q_{2}e^{bv}+r_{2}e^{-bv}, \\ \vspace{1mm} 
Z(w) = p_{3}+q_{3}e^{bw}+r_{3}e^{-bw}, 
\end{array}
\end{eqnarray}
where $p_{i}, q_{i}, r_{i}$, $1 \leq i \leq 3$ are constants. They need to satisfy the equation (6.2) with $w = -u-v$, which gives an equation such as
\[P_{1}e^{-bu}+P_{2}e^{-bv}+P_{3}e^{bu}+P_{4}e^{bv}+P_{5}e^{-bu-bv}+P_{6}e^{bu+bv} = 0. \]
Here the coefficients $P_i$, $1 \leq i \leq 6$ should be zero. Hence, $p_{i}, q_{i}, r_{i}$ need to satisfy the following six equations: 
\begin{eqnarray}
\hspace{5mm} \left\{ \begin{array}{c}
(p_{2}+p_{3})r_{1}-2q_{2}q_{3} = 0, \\
(p_{1}+p_{3})r_{2}-2q_{1}q_{3} = 0, \\
(p_{2}+p_{3})q_{1}-2r_{2}r_{3} = 0, \\
(p_{1}+p_{3})q_{2}-2r_{1}r_{3} = 0, \\
(p_{1}+p_{2})q_{3}-2r_{1}r_{2} = 0, \\
(p_{1}+p_{2})r_{3}-2q_{1}q_{2} = 0. 
\end{array} \right. 
\end{eqnarray}

\vspace{2mm}

(ii) The case where $a < 0$. Let $a = -b^2$ for $b > 0$. The solutions of (6.9) are given by
\begin{eqnarray}
\hspace{5mm} \begin{array}{c}
X(u) = p_{1}+q_{1}\cos{(bu)}+r_{1}\sin{(bu)}, \\ \vspace{1mm} 
Y(v) =p_{2}+q_{2}\cos{(bv)}+r_{2}\sin{(bv)}, \\ \vspace{1mm} 
Z(w) = p_{3}+q_{3}\cos{(bw)}+r_{3}\sin{(bw)}. 
\end{array}
\end{eqnarray}
Substituting them into (6.2) with $w = -u-v$ and using addition formulas, we have an equation such as
\[P_{1}\sin{(bu)}+P_{2}\cos{(bu)}+P_{3}\sin{(bv)}+P_{4}\cos{(bv)} \]
\[+P_{5}\sin{(bu+bv)}+P_{6}\cos{(bu+bv)} = 0. \]
So $P_i = 0$ and $p_{i}, q_{i}, r_{i}$ need to satisfy the following equations: 
\begin{eqnarray}
\hspace{5mm} \left\{ \begin{array}{c}
(p_{2}+p_{3})q_{1}-q_{2}q_{3}+r_{2}r_{3} = 0, \\
(p_{2}+p_{3})r_{1}+q_{2}r_{3}+q_{3}r_{2} = 0, \\
(p_{1}+p_{3})q_{2}-q_{1}q_{3}+r_{1}r_{3} = 0, \\
(p_{1}+p_{3})r_{2}+q_{1}r_{3}+q_{3}r_{1} = 0, \\
(p_{1}+p_{2})q_{3}-q_{1}q_{2}+r_{1}r_{2} = 0, \\
(p_{1}+p_{2})r_{3}+q_{1}r_{2}+q_{2}r_{1} = 0. 
\end{array} \right. 
\end{eqnarray}

\vspace{2mm}

(iii) The case where $a = 0$. The solutions of (6.9) are given by
\begin{eqnarray}
\hspace{5mm} \begin{array}{c}
X(u) = p_{1}+q_{1}u+r_{1}u^2, \\ \vspace{1mm} 
Y(v) =p_{2}+q_{2}v+r_{2}v^2, \\ \vspace{1mm} 
Z(w) = p_{3}+q_{3}w+r_{3}w^2. 
\end{array}
\end{eqnarray}
Substituting them into (6.2) with $w = -u-v$, we have an equation such as
\[P_{1}+P_{2}u+P_{3}v+P_{4}u^{2}+P_{5}v^{2}+2(P_{4}+P_{5})uv+P_{6}(u^{2}v+uv^{2}) = 0. \]
So $P_i = 0$ and $p_{i}, q_{i}, r_{i}$ need to satisfy the following equations: 
\begin{eqnarray}
\hspace{5mm} \left\{ \begin{array}{c}
(p_{2}+p_{3})q_{1}+(p_{1}+p_{3})q_{2}+(p_{1}+p_{2})q_{3} = 0, \\
2(p_{2}+p_{3})r_{1}-2(p_{1}+p_{2})r_{3}+q_{2}(q_{1}-q_{3}) = 0, \\
2(p_{1}+p_{3})r_{2}-2(p_{1}+p_{2})r_{3}+q_{1}(q_{2}-q_{3}) = 0, \\
(q_{2}-q_{3})r_{1}-(q_{1}-q_{2})r_{3} = 0, \\
(q_{1}-q_{3})r_{2}+(q_{1}-q_{2})r_{3} = 0, \\
r_{1}r_{2}+r_{1}r_{3}+r_{2}r_{3} = 0. 
\end{array} \right. 
\end{eqnarray}

\begin{thm}
Let $M$ be a separable minimal surface in $({\mathbb R}^3, \|\cdot\|_{2m})$ as above. 

(i) The case $a > 0$. The functions $X$, $Y$ and $Z$ are given by (6.10) where $p_{i}, q_{i}, r_{i}$ satisfy (6.11).

(ii) The case $a < 0$. The functions $X$, $Y$ and $Z$ are given by (6.12) where $p_{i}, q_{i}, r_{i}$ satisfy (6.13). 

(iii) The case $a = 0$. The functions $X$, $Y$ and $Z$ are given by (6.14) where $p_{i}, q_{i}, r_{i}$ satisfy (6.15). 
\end{thm}

Conversely, from $X$, $Y$ and $Z$, we have
\[x_1 = \pm\int^{u} (X(u))^{-\frac{2m-1}{2m}} du, \ \ \ x_2 = \pm\int^{v} (Y(v))^{-\frac{2m-1}{2m}} dv, \]
\[x_3 = \pm\int^{w=-u-v} (Z(w))^{-\frac{2m-1}{2m}} dw, \]
where for $x_3$, $w = -u-v$ is substituted after integration. 

Now, choosing $p_{i}, q_{i}, r_{i}$ suitably, we will get separable minimal surfaces in $({\mathbb R}^3, \|\cdot\|_{2m})$, where $X$, $Y$ and $Z$ should be positive, and the domain should be non-empty.

\vspace{2mm}

{\bf Example 6.1.} The case $a = 0$. Choose 
\[p_1 = -1, \ \ \ q_1 = 0, \ \ \ r_1 = 1, \]
\[p_2 = -1, \ \ \ q_2 = 0, \ \ \ r_2 = 1, \]
\[p_3 = 2, \ \ \ q_3 = 0, \ \ \ r_3 = -\frac{1}{2}. \]
Then (6.15) holds, and
\[X(u) = u^{2}-1, \ \ \ Y(v) = v^{2}-1, \ \ \ Z(w) = 2-\frac{1}{2}w^{2}. \]
The domain is given by
\[|u| > 1, \ \ \ |v| > 1, \ \ \ |u+v| < 2, \]
and
\[x_1 = \pm\int^{u} (u^{2}-1)^{-\frac{2m-1}{2m}} du, \ \ \ \ x_2 = \pm\int^{v} (v^{2}-1)^{-\frac{2m-1}{2m}} dv, \]
\[x_3 = \pm\int^{w=-u-v} \left( 2-\frac{w^2}{2} \right)^{-\frac{2m-1}{2m}} dw. \]

\vspace{2mm}

{\bf Example 6.2.} The case $a =1$. Choose 
\[p_1 = 1, \ \ \ q_1 = 0, \ \ \ r_1 = 1, \]
\[p_2 = 1, \ \ \ q_2 = 0, \ \ \ r_2 = 1, \]
\[p_3 = -1, \ \ \ q_3 = 1, \ \ \ r_3 = 0. \]
Then (6.11) holds, and
\[X(u) = 1+e^{-u}, \ \ \ Y(v) = 1+e^{-v}, \ \ \ Z(w) = e^{w}-1. \]
The domain is given by $u+v < 0$, and
\[x_1 = \pm\int^{u} (1+e^{-u})^{-\frac{2m-1}{2m}} du, \ \ \ \ x_2 = \pm\int^{v} (1+e^{-v})^{-\frac{2m-1}{2m}} dv, \]
\[x_3 = \pm\int^{w=-u-v} (e^{w}-1)^{-\frac{2m-1}{2m}} dw. \]

\vspace{2mm}

{\bf Example 6.3.} The case $a =-1$. Choose 
\[p_1 = 1, \ \ \ q_1 = 0, \ \ \ r_1 = 1, \]
\[p_2 = 1, \ \ \ q_2 = 0, \ \ \ r_2 = -1, \]
\[p_3 = -\frac{1}{2}, \ \ \ q_3 = \frac{1}{2}, \ \ \ r_3 = 0. \]
Then (6.13) holds, but
\[X(u) = 1+\sin{u}, \ \ \ Y(v) = 1-\sin{v}, \ \ \ Z(w) = -\frac{1}{2}+\frac{1}{2}\cos{w} \leq 0. \]
So a separable minimal surface is not given for such a choice of $p_{i}, q_{i}, r_{i}$.

\vspace{3mm}

To find an example for a negative $a$, we will rewrite (6.13). First, we note that, by (6.12), for the positivity of $X$, $Y$ and $Z$, the following conditions are necessary: 
\[p_{i}+\sqrt{q_{i}^{2}+r_{i}^{2}} > 0, \ \ \ 1 \leq i \leq 3. \]
Set
\[s_{i} := q_{i}+r_{i}\sqrt{-1}, \ \ \ 1 \leq i \leq 3. \]
Then (6.13) is equivalent to that
\begin{eqnarray}
(p_{2}+p_{3})s_{1} = \bar{s_2}\bar{s_3}, \ \ \ (p_{1}+p_{3})s_{2} = \bar{s_1}\bar{s_3}, \ \ \ (p_{1}+p_{2})s_{3} = \bar{s_1}\bar{s_2}. 
\end{eqnarray}
So, if $p_i$, $1 \leq i \leq 3$ are given, then we should have
\[s_1 = \frac{\bar{s_2}\bar{s_3}}{p_{2}+p_{3}}, \ \ \ |s_2|^{2} = (p_{1}+p_{2})(p_{2}+p_{3}), \ \ \ |s_3|^{2} = (p_{1}+p_{3})(p_{2}+p_{3}), \]
and $p_{i}+|s_i| > 0$ for $1 \leq i \leq 3$. 

Noticing the above, we have an example in the case where $a < 0$.

\vspace{2mm}

{\bf Example 6.4.} The case $a =-1$. Let $p_1 = p_2 = p_3 = 1$. Then we should have $|s_2|^2 = |s_3|^2 = 4$. For example, let us choose $s_2 = s_3 = \sqrt{2}+\sqrt{2} \sqrt{-1}$, that is, $q_2 = r_2 = q_3 = r_3 = \sqrt{2}$. Then we should have $s_1 = -2\sqrt{-1}$, that is, $q_1 = 0$ and $r_1 = -2$. For such $p_{i}, q_{i}, r_{i}$, we can check that (6.13) holds, and
\[X(u) = 1-2\sin{u}, \ \ \ Y(v) = 1+\sqrt{2}\cos{v}+\sqrt{2}\sin{v}, \]
\[Z(w) = 1+\sqrt{2}\cos{w}+\sqrt{2}\sin{w}. \]
The domain is given by the conditions
\[\sin{u} < \frac{1}{2}, \ \ \ \sin{\left( v+\frac{\pi}{4} \right)} > -\frac{1}{2}, \ \ \ \sin{\left( \frac{\pi}{4}-u-v \right)} > -\frac{1}{2}, \]
which is non-empty, and
\[x_1 = \pm\int^{u} \left( 1-2\sin{u} \right)^{-\frac{2m-1}{2m}} du, \]
\[x_2 = \pm\int^{v} \left( 1+\sqrt{2}\cos{v}+\sqrt{2}\sin{v} \right)^{-\frac{2m-1}{2m}} dv, \]
\[x_3 = \pm\int^{w=-u-v} \left( 1+\sqrt{2}\cos{w}+\sqrt{2}\sin{w} \right)^{-\frac{2m-1}{2m}} dw. \]

\vspace{2mm}

Graduate School of Science and Technology 

Hirosaki University 

Hirosaki 036-8561, Japan 

E-mail: sakaki@hirosaki-u.ac.jp 

\end{document}